\documentclass[11pt,reqno]{amsart}
\usepackage{amssymb,amscd,amsbsy}
\begin{document}

\newcommand{\vse}{\vspace{.2in}}
\numberwithin{equation}{section}

\title{\bf Utilization of technology for mathematical talks. An alarming situation}
\author{V.V. Peller}

\maketitle

I believe all mathematicians will agree that over the last ten years
the percentage of mathematical talks given by a computer
presentation has jumped and is still rising. I find this situation
alarming. I have the strong opinion that this tendency leads to a
degradation of the quality of talks and discredits the idea of
mathematical talks.

Unfortunately, it is not rare for the organizers of mathematical
conferences to write to the participants to encourage them to give a
computer presentation or, even worse, they inform speakers that
there will be no blackboards at all. On many occasions I have asked
the organizers why they do not provide blackboards. They always
responded that the conference would be held in a hotel (or another
place) that has no blackboards. I do not find this explanation
satisfactory. {\it I strongly believe that if they want to organize
a conference, it is their responsibility to provide blackboards. No
excuse is acceptable.} Incidentally, twenty years ago most
conferences were held at universities, while presently it is more
common to organize conferences at hotels.

Once I was invited to give a seminar talk and I was informed that
they have no blackboards and I should prepare a computer
presentation. I said that was absolutely impossible. It is as if you
invite a pianist to give a performance, but unfortunately you have
no piano and you encourage the pianist to record the notes onto a
flash drive for their computer to be able to reproduce the notes.
After my words a blackboard was found.

\

\begin{center}
{\bf\large Is a mathematician a performer or a composer?}
\end{center}

\

Once giving a talk at a conference, I mentioned the above comparison
with inviting a pianist not having a piano. One mathematician
disagreed with my comparison; he said that a mathematician is a
composer, not a performer.

I strongly disagree with his opinion. When a mathematician proves
theorems, he is a composer. However, when he gives mathematical
talks, {\it he is a performer}. Mathematicians should not
underestimate the importance of being a performer.

\

\begin{center}
{\bf\large What is wrong with computer presentations?}
\end{center}

\

Most speakers who give a computer presentation completely lose
control over the speed of presentation. They switch the screen very
frequently and they completely forget that the audience need time to
digest the information, to comprehend the material. As a result, it
is very difficult (or even impossible) to follow the speaker.

Once I attended a plenary one hour lecture. During the lecture the speaker stated all his results for the last three years and all the results of his students for the last three years. The speaker switched the screen over a hundred times. No need to say that it was absolutely impossible to grasp the picture. Usually, when the speaker loses the audience he does not realize this, he himself watches the screen and does not care about the audience.

Very often the speaker takes the pdf file of his paper designed for publication and displays it on the screen. This is a complete disaster. The speaker understands that it is impossible for the audience to read all the paper, he wants to skip certain slides, he moves the material back and forth by trying to find a page he wants to concentrate on. The effect is disastrous. The only thing the audience can get from such a lecture is blinking of screens.

Sometimes the speaker displays material with a detailed text. At the
same time he talks and gives a verbal explanation of the material
and his words do not coincide with the text on the screen.
Apparently, he believes that people can at the same time read the
text on the screen and listen to another text. Myself, I am not able
to do two such things at the same time.

Even if the speaker cares about the audience, goes at a reasonable
speed and thoroughly prepares the slides, there are still certain
things that create problems for the audience. For example, when the
speaker switches the screen, everything on the previous screen is
gone. However, for people in the audience it is often important to
look again at a formula on the previous screen. Sometimes speakers
combine the computer and the blackboard, and write important things
on the blackboard. This can improve the presentation. However, this
leads to the following problem: to read the screen, they turn off
the light; when the speaker writes on the blackboard, they turn on
the light. This can be very distractive.

Another problem is that when the speaker shows the new screen, everything on the screen
appears at the same time. It is awkward. It is much easier to comprehend the material if the contents of the screen appear gradually, at a speed at which it can be digested. It is true that sometimes computer speakers renew the screen in small portions. This can slightly improve the quality of the presentation. Nevertheless, it is still impossible to display a complicated formula gradually. However, for the audience it is much easier to digest the material if a formula appears on the screen gradually in front of its eyes.

I strongly believe that at least to a certain extent the talk should be an improvisation. The speaker should not lose the audience; he has to have  feedback from the audience. The feedback should help the speaker to adjust the speed of the presentation and decide whether he can include in the talk everything he planned or should omit certain parts of the talk (or, which is very very rare, sometimes to add some material). As Paul Halmos said in \cite{H},
"The faces in the audience can be revealing and
helpful: they can indicate the need to slow down, to speed up, to explain something, to omit something."
It is very difficult to achieve this if you give a computer presentation. Once the file is prepared, the speaker cannot alter it during the talk.

Let me also mention technical problems. From time to time (usually such things happen at least once at each conference) the speaker brings his flash, inserts it in the computer and the computer for some (mysterious) reason does not want to work with his file. Sometimes it takes 10 to 30 minutes to fix the problem.
Another common technical problem: pop-ups appear on the screen, and the speaker desperately attempts to delete them, often asking for assistance.

I do not want to say that all computer presentations are bad. There are exceptions. There are speakers who can minimize the negative features of computer presentations. However, if they were to use their skills to give a talk at the blackboard, the result would be better. I certainly do not want to say that all blackboard talks are good. Not at all! There are speakers who give blackboard talks very poorly. However, I am a hundred percent sure that if those speakers were to give computer presentations, the result would be even worse.

\

\begin{center}
{\bf\large Why do people want to give computer presentations?}
\end{center}

\

I think for many mathematicians the idea of giving a computer
presentation is attractive because they believe that it is easier
for them to prepare a written text and then follow that text during
the talk. Some speakers feel themselves more confident if they
prepare a computer talk and believe that they have already done all
the job. The remaining part is easy: just follow the file prepared
on the computer. This is certainly a wrong perception. The speaker
must think what is better for the audience, not for himself.

Also, many mathematicians give the same talk on several different
occasions. In this case they need prepare the talk only once.

Several mathematicians told me that they believe that computer talks
have an advantage because the fact that the speaker prepared a
computer talk implies that he prepared himself to the talk. I
strongly disagree with this opinion.  On the contrary, it is often
the case that the speaker having prepared a computer talk, believes
that he is fine which is often wrong. Even if the file is prepared,
the most important and difficult task is to speak in front of the
audience.

I have heard from various mathematicians that for plenary lectures
they use the blackboard, while in the case of short communications
they have to give computer presentations, because otherwise they
have no chance of covering all the material. Again, a wrong
perception. They can certainly cover more material in a computer
presentation, but this does not mean that people in the audience
will be able to digest the material. It is unnecessary to cover
everything related to the subject of your talk. You have to select
the material that can be understood by the audience and get the
audience interested in what you are doing. People who become
interested in your work can learn the details later. To attract
people to your research, it is not a good idea to tell them
everything you can for 20 minutes. I am afraid in this way you will
rather repel the audience from your work. In this respect I would
like to quote Paul Halmos \cite{H}: "If someone told you, in half an
hour, the meaning of each ideogram on a page of Chinese, could you
then read and enjoy the poem on that page in the next half hour?"
The speaker should not try to give too much detail. Sometimes
instead of writing complicated formulae or giving sophisticated
definitions, it is more appropriate to describe verbally what one is
doing.

Sometimes speakers try to give detailed proofs of certain results
and believe that they do not have enough time for the proof if they
give a blackboard talk. The same objection applies. It is not
appropriate to give too many details during the talk.

\

\begin{center}
{\bf\large Is the overhead projector a better alternative?}
\end{center}

\

My opinion is that overhead projector talks are {\it slightly} better
than computer presentations. Usually speakers who use overhead projecting don't
violate the speed limit as badly as certain speakers who give
computer presentations. It is much less common among users of
overhead projecting to display their papers designed for
publication. However, with some speakers the speed is still
excessive; when the screen is changed, everything on that screen is
gone forever and the audience has no chance to review the contents
of that screen. When the lecture room has two overhead projectors,
the speaker can keep some important parts of his talk on one screen
and change the contents of the other screen. This can slightly
improve the quality of the presentation.

However, my opinion is that there are still many disadvantages of
overhead projecting. Let me again quote Paul Halmos: "Do not, ever,
greet an audience with a carefully prepared blackboard (or overhead
projector sheets) crammed with formulas, definitions, and theorems.
(An occasionally advisable exception to this rule has to do with
pictures -if a picture, or two pictures, would help your exposition
but would take too long to draw as you talk, at least with the care
it deserves, the audience will forgive you for drawing it before the
talk begins.) The audience can take pleasure in seeing the visual
presentation grow before its eyes - the growth is part of your
lecture, or should be." I entirely agree with his opinion.

\

\begin{center}
{\bf\large Is it good when students teach their professors?}
\end{center}

\

I heard from several colleagues that they never prepare computer
presentations. However, occasionally when they had joint results
with their students, the students gave computer talks and then sent
the computer files to their professors. Later the professors
themselves gave computer talks using the files of their students.
This is ridiculous! Instead of teaching students how to give talks,
professors are taught by their students.

Recently I attended a computer presentation at a conference. After the talk
a student of the speaker, being in the audience, said that the speaker
gave a computer presentation for the first time in his life, he had given only
blackboard talks before. The student was very proud of his professor and 
proposed to congratulate him on this occasion. {\it I do not think that such a congratulation was a good idea}.

It is common for young people to have a desire to master a modern
technology. Sometimes this is good (in particular, it is certainly
good to be able to prepare a mathematical paper with the help of a
computer). However, this is not always the case. For example, young
children should not use calculators to perform arithmetical
operations, while college students should not use sophisticated
calculators to differentiate elementary functions.

\

\begin{center}
{\bf\large Our old friend the blackboard}
\end{center}

\

I strongly believe that the blackboard is the best medium for mathematical talks. It has many advantages over computer presentations or overhead projecting. On the other hand, I do not know any disadvantage (I speak from the point of view of people in the audience; the speaker can say that for him it is more convenient to prepare a computer presentation). Whatever can be performed with the help of a computer or overhead projector, can also be performed much better on the blackboard.

First of all, the speed of a blackboard presentation is limited by
the speed of writing on the board. I do not want to say that it is
impossible to speak too fast. It is. However, not to the same extent
as many speakers do with computer presentations.

Next, the speaker does not have to erase the material as frequently
as computer speakers switch the screen. It is certainly much better
if the lecture room has several big blackboards. In this case the
speaker can erase the contents of a board very rarely and in the
case of short communications there is no need to erase anything from
the boards. Even if the speaker has to use an eraser, he can reserve
one of the blackboards for the material that will not be erased
during the talk.

What else? I have already said that the audience can digest the material much better if the material appears in front of its eyes gradually. It is much more difficult to comprehend the material if a large portion of it is suddenly displayed.

Finally, many computer (or overhead projector) speakers always look
at the screen.  The speaker at the blackboard does not look at the
screen. He can watch the audience instead. As the result he can have
a better contact with the audience, have a feedback from the
audience and be sure that he has not lost the audience. As I have
already said above, having a good contact with the audience allows
the speaker to adjust the speed of presentation and the amount of
material to be presented. All this is very important.

Again, let me repeat that not all blackboard talks are good talks.
The speaker has to know how to use the advantages of the blackboard
to be able to give a good talk. Nevertheless, I have no doubt that
the percentage of good blackboard talks is much higher than the
percentage of good computer (or overhead projecting) talks.

\

\begin{center}
{\bf\large What can be done to rescue conferences?}
\end{center}

\

This is a tough question. My opinion is that only the blackboard should be used for mathematical talks. I know that some mathematicians will call me an extremist. However, even if the mathematical community is not ready to get rid of computer presentations, there are several steps that can be taken.

First of all, I think mathematicians who use computers to give talks should ask themselves why they do it. Perhaps they have no reason other than that it is more convenient for them to prepare a computer presentation. If this is the only reason, they should think about the audience rather than their own convenience. If they do it because otherwise they will not be able to cover all the material, they should understand that using a computer, they can increase the speed of the presentation, but they {\it cannot} increase the speed of comprehension.

I think it would be a good idea for PhD advisors to teach their students how to give talks and explain to them that if they learn how to give good talks at the blackboard, they will produce a much better impression on the audience, which will help their career.

Finally, I believe that the mathematical community ought to agree on a strict rule: {\it organizers of conferences must provide several blackboards to each lecture room}. Also, the organizers of conferences should by no means encourage participants to give computer presentations.

\medskip

{\bf Acknowledgement.} I would like to thank Nicholas Young for reading the manuscript and for making helpful suggestions.

\

\noindent
\begin{tabular}{p{8cm}p{14cm}}
V.V. Peller \\
Department of Mathematics \\
Michigan State University  \\
East Lansing, Michigan 48824\\
USA
\end{tabular}

\end{document}